\documentclass[12pt]{article}
\usepackage[utf8]{inputenc}
\usepackage{mathtools}
\usepackage{booktabs}
\usepackage[english]{babel} % English language/hyphenation% Better typography
\usepackage{amsmath,amsfonts,amsthm}
\usepackage{ amssymb }
\usepackage{graphicx}
\usepackage{array}
\usepackage{multirow}
\usepackage{hhline}
\usepackage[titletoc]{appendix}
\usepackage{subfigure}

\usepackage{booktabs} % Horizontal rules in tables
\usepackage{natbib}
\usepackage{setspace}
\usepackage[margin=1in]{geometry}

\usepackage{subfig}
\usepackage{tabularx}
\usepackage[font = small,labelfont=bf,textfont=it]{caption} % Custom captions under/above floats in tables or figures
\usepackage{footnote}
\usepackage{algorithm}% http://ctan.org/pkg/algorithms
\usepackage[noend]{algpseudocode}% http://ctan.org/pkg/algorithmicx

\newtheorem{theorem}{Theorem}

\newtheorem{lemma}{Lemma}
\newtheorem{definition}{Definition}
\newtheorem{assumption}{Assumption}
\newtheorem{proposition}{Proposition}

\theoremstyle{definition}
\newtheorem{example}{Example}

\newcommand{\RR}{\mathbb{R}}

\newcommand{\vertiii}[1]{||| #1 |||}

\begin{document}
\title {Reduced Rank Multivariate Kernel Ridge Regression}% change the number accordingly.
\author{Wenjia Wang\\The Statistical and Applied Mathematical Sciences Institute\\
Durham, NC, U.S.A. \\wenjia.wang234@duke.edu\\ \and Yi-Hui Zhou\\Department of Biological Sciences\\
North Carolina State University, Raleigh, NC, U.S.A. \\ yihui$\_$zhou@ncsu.edu}
\date{}% do not change this line
\maketitle
\vspace{-5mm}

\begin{abstract}
In the multivariate regression, also referred to as multi-task learning in machine learning, the goal is to recover a vector-valued function based on noisy observations. The vector-valued function is often assumed to be of low rank. Although the multivariate linear regression is extensively studied in the literature, a theoretical study on the multivariate nonlinear regression is lacking. In this paper, we study reduced rank multivariate kernel ridge regression, proposed by \cite{mukherjee2011reduced}. We prove the consistency of the function predictor and provide the convergence rate. An algorithm based on nuclear norm relaxation is proposed. A few numerical examples are presented to show the smaller mean squared prediction error comparing with the elementwise univariate kernel ridge regression.
\end{abstract}

% \begin{keyword}[class=MSC]
% \kwd[Primary ]{62K99}
% \kwd{62M30}
% \kwd[; secondary ]{62M20}
% \end{keyword}

% \begin{keywords}
%   Multivariate regression, Multi-task learning, Reduced rank, Kernel ridge regression
% \end{keywords}

% \end{frontmatter}

\section{Introduction}

Multivariate regression, also referred to as multi-task learning in machine learning, is widely used in machine learning \citep{zhen2017multi}, biology \citep{zhang2012multi}, ecology \citep{kocev2009using}, chemometrics \citep{burnham1999latent}, energy systems \citep{han2012real}, etc. In multivariate regression, the goal is to recover a vector-valued function given the scattered data. Because many popular methods are related to multivariate regression, including principal components \citep{massy1965principal}, partial least squares \citep{wold1975soft}, and vector autoregressive processes  \citep{lutkepohl2013introduction}, it has gained considerable attention from statistics. An incomplete list of works includes \cite{yuan2007dimension,negahban2011estimation,peng2010regularized,frank1993statistical,bedrick1994model,fujikoshi1997modified,lutz2006boosting,mishra2014fixed} and the references therein. In these works, the responses are modeled as linear combinations of the inputs.

In some cases, the responses may have more complicated structure which is nonlinear. Therefore, multivariate linear regression is not adequate to capture the nonlinear structure. To address nonlinear problems, nonlinear multi-task learning methods have been proposed. For example, \cite{bonilla2008multi} propose a multi-task Gaussian prediction, while \cite{liao2006radial} propose a radial basis function neural network. Deep neural networks for multi-task learning are also proposed in the literature, see \cite{zhang2014facial,liu2015multi, zhang2016deep, li2014heterogeneous} for example.

Despite of the wide use of nonlinear multi-task learning, there are very few theoretical studies on multivariate nonlinear regression. Like other machine learning based methods, the nonlinear multi-task learning methods mentioned above are not fully understood from a theoretical point of view for convergence and consistency. To the best of our knowledge, the only theoretical study of multivariate nonlinear regression is \cite{foygel2012nonparametric}, where an additive model is considered, and the goal is to recover a matrix on sample points, instead of recovering the vector-valued function. Therefore, a thorough study of multivariate nonlinear regression is still lacking.

In this work, we study the reduced rank multivariate kernel ridge regression proposed by \cite{mukherjee2011reduced}. This approach is a natural extension of the univariate kernel ridge regression. To the best of our knowledge, the reduced rank multivariate kernel ridge regression, and all other nonparametric multivariate regression, do not have a theoretical guarantee for the estimation of the underlying function. Therefore, our work is the first work providing a theoretical study on the prediction consistency of the nonparametric multivariate regression. We show that, the reduced rank multivariate kernel ridge regression has the same convergence rate as the elementwise univariate kernel ridge regression when the underlying function is nonsparse, and has a faster convergence rate when the underlying function is sparse. We propose an algorithm for solving the optimization problem in the reduced rank multivariate kernel ridge regression by nuclear norm relaxation.

\section{Background and methodology}\label{secBackground}
In this section, we provide an introduction to the univariate kernel ridge regression, and present the reduced rank multivariate kernel ridge regression.

\subsection{Univariate kernel ridge regression}
In the nonparametric univariate regression problem, the goal is to recover a univariate function with observed data. Let $\Omega\subset \RR^d$ be the region of interest, which is compact and convex. Suppose we observe data $(x_k,y_k),k=1,\ldots,n$, given by
\begin{eqnarray*}
y_k=f(x_k) + \epsilon_k, %\label{recovering}
\end{eqnarray*}
where $x_k \in \Omega$, $y_k\in \RR$, $f:\Omega \rightarrow \RR$ is a univariate function, and $\epsilon_k$'s are independent and identically distributed random variables with mean zero and finite variance. One widely used method is the univariate kernel ridge regression \citep{saunders1998ridge}, which employs the reproducing kernel Hilbert space $\mathcal{N}_\Psi(\Omega)$ generated by kernel function $\Psi$. The reproducing kernel Hilbert space can be defined via Fourier transform, defined by
$$\mathcal{F}(f)(\omega)=(2\pi)^{-d/2}\int_{\mathbb{R}^d} f(x) e^{-ix^T\omega}d x,$$ for $f\in L_1(\mathbb{R}^d)$. The definition of the reproducing kernel Hilbert space can be generalized to $f\in L_2(\RR^d)\cap C(\RR^d)$. See \cite{girosi1995regularization} and Theorem 10.12 of \cite{wendland2004scattered}.
\begin{definition}\label{Def:NativeSpace}
Let $\Psi$ be a kernel function. Define the reproducing kernel Hilbert space $\mathcal{N}_\Psi(\RR^d)$ generated by $\Psi$ as
	$$\mathcal{N}_\Psi(\RR^d):=\{f\in L_2(\RR^d)\cap C(\RR^d):\mathcal{F}(f)/\sqrt{\mathcal{F}(\Psi)}\in L_2(\RR^d)\},$$
	with the inner product
	$$\langle f,g\rangle_{\mathcal{N}_\Psi(\RR^d)}=(2\pi)^{-d}\int_{\RR^d}\frac{\mathcal{F}(f)(\omega)\overline{\mathcal{F}(g)(\omega)}}{\mathcal{F}(\Psi)(\omega)}d \omega.$$
\end{definition}
A reproducing kernel Hilbert space on $\Omega\subset \RR^d$ is denoted by $\mathcal{N}_\Psi(\Omega)$, with norm
\begin{eqnarray*}%\label{restriction}
\|f\|_{\mathcal{N}_\Psi(\Omega)}=\inf\{\|f_E\|_{\mathcal{N}_\Psi(\RR^d)}:f_E\in\mathcal{N}_\Psi(\RR^d),f_E|_\Omega=f\},
\end{eqnarray*}
where $f_E|_\Omega$ denotes the restriction of $f_E$ to $\Omega$.

A prominent class of kernel functions is the isotropic Mat\'ern kernel function \citep{stein2012interpolation}, which is given by
\begin{align}\label{matern}
	\Psi(x)=\frac{1}{\Gamma(\nu - d/2)2^{\nu - d/2 -1}}\|x\|_2^{\nu - d/2} K_{\nu - d/2}(\|x\|_2)
\end{align}
after a proper reparametrization, where $K_{\nu - d/2}$ is the modified Bessel function of the second kind, and $\|\cdot\|_2$ denotes the Euclidean metric. The parameter $\nu>d/2$ is the smoothness parameter, because it can control the smoothness of functions in $\mathcal{N}_{\Psi}(\Omega)$. It can be shown that the reproducing kernel Hilbert space generated by the isotropic Mat\'ern kernel function coincides with the Sobolev space with smoothness $\nu$, denoted by $H^\nu(\Omega)$. For a positive number $\nu > d/2$, the Sobolev space on $\RR^d$ with smoothness $\nu$ can be defined as
\begin{align*}
H^\nu(\mathbb{R}^d) = \{f\in L_2(\mathbb{R}^d): |\mathcal{F}(f)(\cdot)| (1+\|\cdot\|^2)^{\nu/2}\in L_2(\mathbb{R}^d)\},
\end{align*}
equipped with an inner product
$$\langle f,g\rangle_{H^\nu(\mathbb{R}^d)}=(2\pi)^{-d}\int_{\RR^d}\mathcal{F}(f)(\omega)\overline{\mathcal{F}(g)(\omega)}(1+\|\omega\|^2)^{\nu}d \omega.$$ The following lemma, which is a direct result of Corollary 10.13 in \cite{wendland2004scattered} and the extension theorem \citep{devore1993besov}, states
that the reproducing kernel Hilbert space $\mathcal{N}_{\Psi}(\Omega)$ coincides with the Sobolev space with smoothness $\nu$ $H^{\nu}(\Omega)$.
\begin{lemma}\label{coro1013}
	Let $\Psi$ be as in \eqref{matern}. We have the following.
	\begin{enumerate}
	    \item The reproducing kernel Hilbert space $\mathcal{N}_{\Psi}(\RR^d)$ coincides with the Sobolev space with smoothness $\nu$ $H^{\nu}(\RR^d)$, and the norms $\|\cdot\|_{\mathcal{N}_{\Psi}(\RR^d)}$ and $\|\cdot\|_{H^\nu(\RR^d)}$ are equivalent.
	    \item Suppose $\Omega$ is compact and convex. Then the reproducing kernel Hilbert space $\mathcal{N}_{\Psi}(\Omega)$ coincides with the Sobolev space with smoothness $\nu$ $H^{\nu}(\Omega)$, and the norms $\|\cdot\|_{\mathcal{N}_{\Psi}(\Omega)}$ and $\|\cdot\|_{H^\nu(\Omega)}$ are equivalent.
	\end{enumerate}
\end{lemma}
The univariate kernel ridge regression reconstructs $f$ by using
\begin{align}\label{KRRestuni}
\hat f = \operatorname*{argmin}_{\hat f\in \mathcal{N}_{\Psi}(\Omega)}\bigg( \frac{1}{n}\sum_{j=1}^n(y_j-\hat f(x_j))^2 + \lambda \|\hat f\|^2_{\mathcal{N}_{\Psi}(\Omega)}\bigg),
\end{align}
where $\lambda>0$ is a parameter. It is known that the optimal convergence rate of of the $L_2$ prediction error for the nonparametric regression is $n^{-\nu/(2\nu + d)}$ \citep{stone1982optimal}. The following theorem is a direct result of Theorem 10.2 of \cite{geer2000empirical} and Lemma \ref{lemmaratiosingle} in the Supplementary materials, which states that by choosing appropriate parameter $\lambda$, the univariate kernel ridge regression can achieve the optimal convergence rate. In the rest of this work, we will use the following notation. For two positive sequences $a_n$ and $b_n$, we write $a_n\asymp b_n$ if, for some constants $C,C'>0$, $C\leq a_n/b_n \leq C'$.
% Similarly, we write $a_n\gtrsim b_n$ if $a_n\geq Cb_n$ for some constant $C>0$. For notational simplicity, we will use $C,C',C_1,C_2,...$ and $\eta,\eta_0,\eta_1,...$ to denote the constants, of which the values can change from line to line.
\begin{theorem}\label{thmunikrr}
Let $\lambda \asymp n^{-\frac{2\nu}{2\nu + d}}$. If $f \neq 0$, we have
\begin{align}\label{thmunikrreq1}
    \|f - \hat f\|_{L_2(\Omega)}^2 =O_{P}\left(n^{-\frac{2\nu}{2\nu + d}}\right)\|f\|^2_{\mathcal{N}_{\Psi}(\Omega)},
\end{align}
otherwise we have
\begin{align}\label{thmunikrreq2}
    \|f - \hat f\|_{L_2(\Omega)}^2 =O_{P}\left(n^{-\frac{2\nu}{2\nu + d}}\right).
\end{align}
\end{theorem}
By the representer theorem \citep{gu2013smoothing}, \eqref{KRRestuni} has the following closed form
\begin{eqnarray}\label{eqrecoveringfX}
\hat{f}(x)=\sum_{j=1}^n u_{j}\Psi(x-x_j),
\end{eqnarray}
with $u=(u_{1},\ldots,u_{n})^T$ satisfying $u = (\Psi(X-X)+n\lambda I_n)^{-1}Y$, where $\Psi(X-X) = (\Psi(x_j - x_k))_{jk}$, $Y=(y_1,\ldots,y_n)^T$, and $I_n$ is the identity matrix.

\subsection{Reduced rank multivariate kernel ridge regression}
In the multivariate regression, the goal is to recover a vector-valued function $F = (f_1,...,f_p)^T$, where $f_k:\Omega \rightarrow \RR$, $k=1,...,p$ is defined on a compact and convex region $\Omega \subset\RR^d$. Define $\mathcal{F}={\rm span}\{f_1,...,f_p\}$ be the linear space spanned by $f_1,...,f_p$. In the reduced rank multivariate regression, it is often assumed that the the dimension of $\mathcal{F}$, denoted by dim$(\mathcal{F})$, is small. Let $r = $dim$(\mathcal{F})$. Nevertheless, in this work we do not require $r$ to be small, that is, we can have $1\leq r\leq p$. Suppose we observe data $(x_j,Y_j'),j=1,\ldots,n$ for $x_j \in \Omega$, satisfying
\begin{eqnarray}\label{recovering}
Y_j'=F(x_j) + \epsilon_j,
\end{eqnarray}
where $Y_j' = (y_{j1},...,y_{jp})^T\in \RR^p$, $F(x_j) = (f_1(x_j),...,f_p(x_j))^T$, and $\epsilon_j = (\epsilon_{j1},...,\epsilon_{jp})^T\in \RR^p$. The random noises $\epsilon_{jk}$ are independent and identically distributed, and have mean zero and finite variance. If $f_k(x) = a_k^Tx$ for all $x\in \Omega$, where $a_k \in \RR^d$, then \eqref{recovering} reduces to a multivariate linear regression.

In order to estimate the function $F$, one naive way is to apply elementwise univariate kernel ridge regression to function $f_k$, i.e., solving $p$ optimization problems
\begin{align}\label{unikrrmodel}
    \hat f_k = \operatorname*{argmin}_{\hat f\in \mathcal{N}_{\Psi}(\Omega)}\bigg( \frac{1}{n}\sum_{j=1}^n(y_{jk}-\hat f(x_j))^2 + \lambda \|\hat f\|^2_{\mathcal{N}_{\Psi}(\Omega)}\bigg)
\end{align}
for $k=1,...,p$, where $\lambda>0$ is a tuning parameter. However, this approach may lose the underlying potentially low rank structure of $F$. Therefore, a regression method which utilizes the underlying low rank structure is desired.

In this work, we consider the reduced rank multivariate kernel ridge regression \citep{mukherjee2011reduced}, which is a natural extension of the univariate kernel ridge regression. In the reduced rank multivariate kernel ridge regression, the dimension of the solution is restricted. Specifically, consider imposing a hard threshold on the rank of $\mathcal{F}$, i.e., dim$(\mathcal{F}) \leq r_1$, where $r_1 \geq r$ is the reduced rank constraint. The rank $r_1$ controls the rank of the solution, thus can control the complexity of the predictor. The reduced rank multivariate kernel ridge regression is to find a function $\hat F = (\hat f_1,...,\hat f_p)^T$ which is the solution to the optimization problem
\begin{align}\label{firstmodelXX}
    & \min_{g_1,...,g_p}\frac{1}{pn}\sum_{j=1}^p \|Y_j - g_j(X)\|_2^2 + \lambda_1 \sum_{k=1}^p\|g_k\|_{\mathcal{N}_{\Psi}(\Omega)}^2,\nonumber\\
    \mbox{ s.t. }& {\rm dim}(\mathcal{G}) \leq r_1,
\end{align}
where $\lambda_1>0$ is a parameter, $Y_j = (y_{1j},...,y_{nj})^T$, $g_j(X) = (g_j(x_1),...,g_j(x_n))^T$ for $j=1,...,p$, and $\mathcal{G}={\rm span}\{g_1,...,g_p\}$. The tuning parameters $\lambda_1$ and $r_1$ can be chosen by generalized cross validation (GCV); see \cite{wahba1990spline}. We will discuss solving \eqref{firstmodelXX} in Section \ref{secNuclear}.

\section{Theoretical properties}\label{secTheory}
In this section, we present our main theoretical results of the reduced rank multivariate kernel ridge regression: the $L_2$ prediction consistency of the reduced rank multivariate kernel ridge regression.

A technical assumption is that the errors are sub-Gaussian \citep{geer2000empirical}.
\begin{assumption}\label{ErrorC1}
Suppose $\epsilon_{jk}$'s in \eqref{recovering} are independent and identically distributed random variables satisfying
\begin{align*}%\label{noiseCond}
C_1^2 \{E( e^{|\epsilon_{jk}|^2/C_1^2})-1\}\leq \sigma_0^2
\end{align*}
for some constants $C_1$ and $\sigma_0$.
\end{assumption}
Before presenting the main theorem, we need to introduce some notation. Consider the Cartesian product of $\mathcal{N}_{\Psi}(\Omega)$, denoted by $\mathcal{N}_{\Psi}^{\otimes p}(\Omega)$, given by
\begin{align*}
    \mathcal{N}_{\Psi}^{\otimes p}(\Omega) = \{F = (f_1,...,f_p)^T: f_j \in \mathcal{N}_{\Psi}(\Omega)\}.
\end{align*}
We can equip $\mathcal{N}_{\Psi}^{\otimes p}(\Omega)$ with a norm $\|F\|$ for $F\in \mathcal{N}_{\Psi}^{\otimes p}(\Omega)$, defined by
\begin{align*}
    \|F\|^2 = \sum_{j=1}^p \|f_j\|_{\mathcal{N}_{\Psi}(\Omega)}^2.
\end{align*}
Similarly, we can define the $L_2$ norm of $F$. With an abuse of notation, we use $\|F\|_{L_2(\Omega)}$ to denote the $L_2$ norm of $F\in \mathcal{N}_{\Psi}^{\otimes p}(\Omega)$, defined by
\begin{align*}
    \|F\|_{L_2(\Omega)}^2 = \sum_{j=1}^p \|f_j\|_{L_2(\Omega)}^2.
\end{align*}
In the following theorem, we show the $L_2$ prediction consistency of the estimator $\hat F$ in \eqref{firstmodelXX}.
\begin{theorem}\label{mainthm1}
Suppose Assumption \ref{ErrorC1} is satisfied, and $r_1 \geq r$, where $r_1$ is as in \eqref{firstmodelXX}. Assume $r_1^{\frac{d}{2\nu}} \|F\|^{\frac{2d}{2\nu+d}}n^{\frac{d}{2\nu+d}}p^{-\frac{2d+2\nu}{2\nu+d}}$ goes to infinity. Let
\begin{align*}
    \lambda_1 \asymp n^{-\frac{2\nu}{2\nu + d}}p^{-\frac{d}{2\nu + d}} \|F\|^{-\frac{4\nu}{2\nu+d}}, \quad  S = \max \{p^{\frac{2\nu}{2\nu+d}}\|F\|^{\frac{2d}{2\nu+d}}, r_1\},
\end{align*}
and $\hat F$ be the solution to \eqref{firstmodelXX}. Then we have
\begin{align}\label{thm1rate}
    \|F - \hat F\|_{L_2(\Omega)}^2 =O_{P}\left(S n^{-\frac{2\nu}{2\nu+d}}\right).
\end{align}
Potential dependency of $r$, $r_1$, $\|F\|$, $S$, and $p$ on $n$ is suppressed for notational simplicity.
\end{theorem}
By Theorem \ref{mainthm1}, it can be seen that if $S n^{-2\nu/(2\nu+d)}$ converges to zero, then the reduced rank multivariate kernel ridge regression can provide a consistent predictor. This indicates that in order to get a consistent estimator, the dimension of the responses cannot be very large compared to the sample size. This is similar to the case of the multivariate linear regression, because in the multivariate linear regression, the dimension of responses should also not be large \citep{negahban2011estimation}.

Using the results in Theorem \ref{mainthm1}, we can obtain the convergence rates of predictor $\hat F$ in different cases. We demonstrate these results in Examples \ref{eg1}-\ref{eg3}.

\begin{example}\label{eg1}
Suppose there exist constants $C_1,C_2>0$ such that for all $f_j\in F$, $C_1\leq \|f_j\|_{\mathcal{N}_{\Psi}(\Omega)} \leq C_2$ for $j=1,...,p$. Therefore, $\|F\|^2 \asymp p$, and obviously we have $r_1\leq p$. By Theorem \ref{mainthm1}, we can choose $$\lambda_1 \asymp n^{-\frac{2\nu}{2\nu + d}}p^{-1}$$ to obtain
\begin{align}\label{eg1rate}
    \|F - \hat F\|_{L_2(\Omega)}^2 =O_{P}\left(pn^{-\frac{2\nu}{2\nu + d}}\right).
\end{align}
If we use the elementwise univariate kernel ridge regression, by Theorem \ref{thmunikrr}, it can be seen that the convergence rate is also $pn^{-2\nu/(2\nu+d)}$. This implies the interesting result that the convergence rate does not suffer by our having introduced the low rank structure, and it is desired to use a low rank structure in practice.
\end{example}

\begin{example}\label{eg2}
If $p=1$, the convergence rate in \eqref{thm1rate} reduces to
\begin{align*}%\label{eg2rate}
    \|F - \hat F\|_{L_2(\Omega)}^2 =O_{P}\left(n^{-\frac{2\nu}{2\nu + d}}\right),
\end{align*}
which coincides with the convergence rate of the univariate kernel ridge regression as in Theorem \ref{thmunikrr}. Therefore, the results in Theorem \ref{mainthm1} is a generalization of the results in the univariate kernel ridge regression.
\end{example}

\begin{example}\label{eg3}
As in Theorem \ref{thmunikrr}, it is allowed that $f_j=0$ for some $j\in \{1,...,p\}$ in Theorem \ref{mainthm1}. If $F$ is sparse, i.e., card$(\{j:f_j\neq 0, f_j\in F\}) = s$, where card$(A)$ is the cardinality of the set $A$, and $s = o(p)$. Suppose there exist constants $C_1,C_2>0$ such that for all $k\in \{j:f_j\neq 0\}$, $C_1\leq \|f_k\|_{\mathcal{N}_{\Psi}(\Omega)} \leq C_2$. Therefore, we have $\|F\|^2\asymp s$. If we choose
\begin{align*}
    \lambda_1 \asymp n^{-\frac{2\nu}{2\nu + d}}p^{-\frac{d}{2\nu + d}} s^{-\frac{2\nu}{2\nu+d}},\quad r_1\asymp s,
\end{align*}
by Theorem \ref{mainthm1}, we obtain the convergence rate
\begin{align*}%\label{eg2rate}
    \|F - \hat F\|_{L_2(\Omega)}^2 =O_{P}\left(p^{\frac{2\nu}{2\nu + d}}n^{-\frac{2\nu}{2\nu + d}}s^{\frac{d}{2\nu + d}}\right).
\end{align*}
If we apply the elementwise univariate kernel ridge regression, by Theorem \ref{thmunikrr}, the convergence rate is $pn^{-2\nu/(2\nu+d)}$. Because $s = o(p)$, the convergence rate of the reduced rank multivariate kernel ridge regression is faster than the convergence rate of the elementwise univariate kernel ridge regression.
\end{example}

From Examples \ref{eg1}-\ref{eg3}, we can see that the convergence rate of the reduced rank multivariate kernel ridge regression is not slower than the convergence rate of the elementwise univariate kernel ridge regression. Because of the low complexity afforded by the low rank structure, it is desired to use the reduced rank multivariate kernel ridge regression.

% An important implication of Theorem \ref{mainthm1} is that in the reduced rank multivariate kernel ridge regression, the dimension of the responses can be large, that is, it can be comparable to the sample size.

\section{Nuclear norm relaxation}\label{secNuclear}
For fixed $\lambda_1$ and $r_1$, it can be shown that the optimization problem \eqref{firstmodelXX} has the form \citep{mukherjee2011reduced}
\begin{align*}
    \hat F(x) = \Psi(x - X)(\Psi(X - X) + n\lambda I_n)^{-1}YP_{r_1}
\end{align*}
for any point $x\in \Omega$, where $\Psi(x - X) = (\Psi(x - x_1),...,\Psi(x - x_n))$, $\Psi(X - X) = (\Psi(x_j - x_k))_{jk}$, $I_n$ is the identity matrix, and $Y = (Y_1,...,Y_p)$. The matrix $P_{r_1}\in \RR^{p \times p}$ is the projection matrix to the space spanned by $r_1$ principal eigenvectors of $P = Y^T\Psi(X - X)(\Psi(X - X) + n\lambda I_n)^{-1}Y$. Because we do not know the information of $F$, we cannot always expect that we can have the optimal choice of $\lambda$ and $r_1$ as in Theorem \ref{mainthm1}. Therefore, we can choose $\lambda_1$ and $r_1$ by generalized cross validation (GCV); see \cite{wahba1990spline}. However, solving the optimization problem in the GCV may be difficult, because the variable $r_1$ in the GCV is not continuous. In this work we propose using nuclear norm relaxation, which is used in matrix completion \citep{candes2009exact}.

By Proposition 1 in \cite{mukherjee2011reduced}, which is an application of the representer theorem, the solution to \eqref{firstmodelXX} has the form
\begin{align}\label{repreform}
    \hat f_k(x) = \sum_{j=1}^n\hat \alpha_{kj} \Psi(x - x_j)
\end{align}
for $x\in \Omega$ and $1\leq k \leq p$, where $\hat\alpha_{kj} \in \RR$. Let $\hat A=(\hat \alpha_{kj})_{kj}\in \RR^{p\times n}$. Then by Proposition 2 of \cite{mukherjee2011reduced}, ${\rm dim}(\hat{\mathcal{F}}) \leq r_1$ is equivalent to rank$(\hat A) \leq r_1$. Therefore, we can write the optimization problem \eqref{firstmodelXX} as
\begin{align}\label{model2}
    & \min_{A\in \RR^{p\times n}}\frac{1}{pn}\sum_{j=1}^p \|Y_j - \Psi(X - X)\alpha_j\|_2^2 + \lambda_1 \sum_{j=1}^p\alpha_j^T\Psi(X - X)\alpha_j,\nonumber\\
    \mbox{ s.t. }& {\rm rank}(A) \leq r_1,
\end{align}
where $A = (\alpha_1,...,\alpha_p)^T$ with $\alpha_j = (\alpha_{j1},...,\alpha_{jn})^T$.

The first relaxation is via Lagrange multiplier. Applying Lagrange multiplier to \eqref{model2}, for a tuning parameter $\lambda_2>0$, up to a diffrence of a constant related to $r_1$ and $\lambda_2$, we have a relaxed optimization problem
\begin{align}\label{model3}
    & \min_{A\in \RR^{p\times n}}\frac{1}{pn}\sum_{j=1}^p \|Y_j - \Psi(X - X)\alpha_j\|_2^2 + \lambda_1 \sum_{j=1}^p\alpha_j^T\Psi(X - X)\alpha_j + \lambda_2 {\rm rank}(A).
\end{align}
One widely used approach to relax the penalty term $\lambda_2 {\rm rank}(A)$ is by replacing the rank with the nuclear norm \citep{candes2009exact}. For a matrix $B\in \RR^{m_1\times m_2}$, denote the ordered singular values of $B$ by $\sigma_1(B) \geq \sigma_2(B) \geq ... \geq\sigma_{\min(m_1,m_2)}(B) \geq 0$. The nuclear norm of $B$ is defined by
\begin{align*}
    \vertiii B = \sum_{j=1}^{\min(m_1,m_2)} \sigma_j(B).
\end{align*}
Note that for any full rank matrix $Q\in \RR^{n\times n}$, the rank of $A$ is eqaul to the rank of $AQ$. Thus, \eqref{model3} is equivalent to the following optimization problem
\begin{align*}
    \min_{A\in \RR^{p\times n}}\frac{1}{pn}\sum_{j=1}^p \|Y_j - \Psi(X - X)\alpha_j\|_2^2 + \lambda_1 \sum_{j=1}^p\alpha_j^T\Psi(X - X)\alpha_j + \lambda_2{\rm rank}(AQ).
\end{align*}
Now we relax the penalty term ${\rm rank}(AQ)$ to $\vertiii {AQ}$. Therefore, our final model becomes
\begin{align}\label{finalmodel}
    \min_{A\in \RR^{p\times n}}\frac{1}{pn}\sum_{j=1}^p \|Y_j - \Psi(X - X)\alpha_j\|_2^2 + \lambda_1 \sum_{j=1}^p\alpha_j^T\Psi(X - X)\alpha_j + \lambda_2\vertiii {AQ}.
\end{align}
The algorithm for solving \eqref{finalmodel} is well-studied in the literature, for example, see \cite{mishra2013low}, and we can use GCV to choose $\lambda_1$ and $\lambda_2$. We recommend choosing $Q = \Psi(X - X)$. Note in \eqref{finalmodel}, $\lambda_1$ and $\lambda_2$ are continuous.

\section{Numeric examples}\label{secNum}
In this section we conduct experiments. Let $a_k = 0.5k$ for $k=1,...,r$. We consider  functions $h_k:[0,1]^d \rightarrow \RR$ for $d=1,2$ and $k=1,...,r$, defined by \citep{sun2014balancing}
\begin{align*}
   h_k(x) = \frac{2}{\sqrt{\sum_{j=1}^d(x_j' - 0.1a_k)^2 } + 1} + \frac{0.5}{\sqrt{\sum_{j=1}^d(x_j' - 0.05a_k)^2} + 1},x_j'\in [0,1],
\end{align*}
for $x = (x_1',...,x_d')^T\in [0,1]^d$. The test functions are chosen to be $F = AH$, where $A\in \RR^{p\times r}$ is a matrix and will be specified later, and $H=(h_1,...,h_r)^T$. Suppose we observe data $(x_j,Y_j'),j=1,\ldots,n$ for $x_j \in [0,1]^d$, given by
\begin{eqnarray*}
Y_j'=F(x_j) + \epsilon_j,
\end{eqnarray*}
where $Y_j' = (y_{j1},...,y_{jp})^T\in \RR^p$, and $\epsilon_j = (\epsilon_{j1},...,\epsilon_{jp})^T\in \RR^p$. The random noise $\epsilon_{jk}\sim N(0,0.1^2)$ are independent and identically distributed, where $N(0,\sigma^2)$ is a normal distribution with mean zero and variance $\sigma^2$. The input points $x_j$ are uniformly sampled on the space $[0,1]^d$. We use
\begin{align*}
    \frac{1}{N}\sum_{j=1}^N \|F(x_j'') - \hat F(x_j'')\|_2^2
\end{align*}
to approximate the $L_2$ prediction error $\|F - \hat F\|_{L_2(\Omega)}^2$, where $N=200$ and $x_1'',...,x_N''$ are first 200 points of the Halton sequence \citep{niederreiter1992random}. We use the isotropic Mat\'ern kernel function as in \eqref{matern} with $\nu = 3.5 + d/2$. In this section, we compare the elementwise univariate kernel ridge regression and reduced rank multivariate kernel ridge regression. In both methods, we choose the tuning parameter using independently generated validation data sets of the same sample size as the input points, as in \cite{mukherjee2011reduced}. For the optimization problem \eqref{finalmodel}, we use different starting points to obtain the optimal point.

\begin{example}[Non-sparse functions]\label{simuegns}
We first consider that $F$ is not sparse. Recall that we choose $F = AH$. Let $A = (I_r,B)^T$, where $I_r\in \RR^{r\times r}$ is an identity matrix and $B\in \RR^{r \times (p-r)}$. Let $B=(b_{jk})_{jk}$, where $b_{jk}$'s are uniformly sampled from $[0,1]$. In the simulation studies, we choose different set of $(d,r,p,n)$. For each set of $(d,r,p,n)$, we randomly generate the matrix $B$, and then use elementwise univariate kernel ridge regression and reduced rank multivariate kernel ridge regression to predict the underlying function $F$. The results are shown in Table \ref{tab:egnonsp}.

\begin{table}[h!]
\centering
\begin{tabular}{cccc}
%\\
 \\
\hline
$(d,r,p,n)$ & EUKRR & RRMKRR & Difference\\
\hline
(1, 2, 10, 20)  & 0.0165 & 0.0147 & 0.0018\\
(1, 2, 10,  60)   & 0.0053 & 0.0046 & 0.0007\\
 (1, 4, 20,  60)  & 0.0136 & 0.0136 & 0.00001\\
 (2, 2, 10,  100)  & 0.2875 & 0.2549 & 0.0328\\%zheli
 (2, 4, 30,  100)  & 1.4444 & 1.2365 & 0.2079\\
\hline
\end{tabular}
\label{tab:egnonsp}
\caption{
The approximated $L_2$ prediction error of the elementwise univariate kernel ridge regression (EUKRR) and reduced rank multivariate kernel ridge regression (RRMKRR). The fourth column denotes the difference of the prediction errors, computed by (second column - third column).}
\end{table}
\end{example}

\begin{example}[Sparse functions]\label{simuegs}
We consider that $F$ is sparse. Recall that we choose $F = AH$. Let $A = (I_r,B,B_1)^T$, where $I_r\in \RR^{r\times r}$ is the identity matrix, and $B_1\in \RR^{r\times (p-s)}$ is a zero matrix, and $s\geq r$ is used to control the sparsity of matrix $A$ and $F$. Let $B=(b_{jk})_{jk} \in \RR^{r\times (s - r)}$, where $b_{jk}$'s are uniformly sampled from $[0,1]$. In the simulation studies, we choose different set of $(d,r,s,p,n)$. For each set of $(d,r,s,p,n)$, we randomly generate the matrix $B$, and then use elementwise univariate kernel ridge regression and reduced rank multivariate kernel ridge regression to predict the underlying function $F$. The results are shown in Table \ref{tab:egsp}.

\begin{table}[h!]
\centering
\begin{tabular}{cccc}
%\\
 \\
\hline
$(d,r,s,p,n)$ & EUKRR & RRMKRR & Difference\\
\hline
(1, 2, 4, 10, 20)  & 0.0276 & 0.0130 & 0.0145\\
(1, 2, 4, 10, 60)   & 0.0052 & 0.0045 & 0.0007\\
 (1, 4,  9, 20,  60) & 0.0128 & 0.0123 & 0.0005\\
 (2, 2,  4, 10,  100) & 0.1783 & 0.1418 & 0.0364\\
 (2, 4,  7, 30,  100)  & 0.6207 & 0.5274 & 0.0933\\
\hline
\end{tabular}
\label{tab:egsp}
\caption{
The approximated $L_2$ prediction error of the elementwise univariate kernel ridge regression (EUKRR) and reduced rank multivariate kernel ridge regression (RRMKRR). The fourth column denotes the difference of the prediction errors, computed by (second column - third column).}
\end{table}
\end{example}

From Examples \ref{simuegns} and \ref{simuegs}, we can see that if the underlying function is of low rank, using reduced rank multivariate kernel ridge regression can obtain a smaller $L_2$ prediction error, compared with the elementwise univariate kernel ridge regression. It is not surprising that if the underlying function is sparse, the reduced rank multivariate kernel ridge regression performs better than the elementwise univariate kernel ridge regression, which corroborates the results of Theorem \ref{mainthm1}. However, even if the underlying function is not sparse, using reduced rank multivariate kernel ridge regression still has a smaller $L_2$ prediction error, although it is not obvious under some cases. Therefore, we conjecture that it is always desirable to use the reduced rank multivariate kernel ridge regression, if the underlying function is of low rank.

\section{Discussion}\label{secDiscuss}
In this work we study the reduced rank multivariate kernel ridge regression, which is used to recover a vector-valued function, from a theoretical perspective. Specifically, we obtain the convergence rate of prediction error using the reduced rank multivariate kernel ridge regression, which is the first result of this kind, as far as we know. In this work, we only consider the original optimization problem \eqref{firstmodelXX} and do not prove the consistency of the predictor obtained by the relaxed model \eqref{finalmodel}. The consistency of the later estimator will be pursued in future work.

\appendix
\renewcommand{\thesection}{\Alph{section}}

\numberwithin{equation}{section}
\numberwithin{lemma}{section}
\numberwithin{assumption}{section}
\numberwithin{proposition}{section}

\section{Multivariate reproducing kernel Hilbert space}\label{subsecmrkhs}
In this section, we introduce a function class and its properties, which are used in our theoretical development.
% In Section \ref{subsecmrkhs}, we first introduce the space that is used to study the function class $\mathcal{F} = {\rm span}\{f_1,...,f_p\}$ with dim$(\mathcal{F}) \leq r$, and is used to prove Theorem \ref{mainthm1}.
Consider the Cartesian product of $\mathcal{N}_{\Psi}(\Omega)$, denoted by $\mathcal{N}_{\Psi}^{\otimes p}(\Omega)$, given by
\begin{align*}
    \mathcal{N}_{\Psi}^{\otimes p}(\Omega) = \{F = (f_1,...,f_p)^T: f_j \in \mathcal{N}_{\Psi}(\Omega),j=1,...,p\}.
\end{align*}
We can equip $\mathcal{N}_{\Psi}^{\otimes p}(\Omega)$ with a norm $\|F\|$ defined by
\begin{align}
    \|F\|^2 = \sum_{j=1}^p \|f_j\|_{\mathcal{N}_{\Psi}(\Omega)}^2,
\end{align}
for any $F\in \mathcal{N}_{\Psi}^{\otimes p}(\Omega)$. Let $\mathcal{A}_r$ be a function class defined by
\begin{align}
    \mathcal{A}_r = \{F = (f_1,...,f_p)^T: F\in \mathcal{N}_{\Psi}^{\otimes p}(\Omega), {\rm dim(span}\{f_1,...,f_p\}) \leq r\}
\end{align}
for $r=1,...,p$. Obviously, $\mathcal{A}_r \subset \mathcal{N}_{\Psi}^{\otimes p}(\Omega)$. Note that $\mathcal{A}_r$ is not a linear subspace, because for $F,F'\in \mathcal{A}_r$, $F + F'$ may not have rank $r$, thus may not be in $\mathcal{A}_r$.

For $R>0$, let
\begin{align*}
    \mathcal{A}_r(R) = \{F\in \mathcal{A}_r: \|F\|\leq R\}.
\end{align*}
We define
\begin{align*}
    \mathcal{B}_r = \{F = (f_1,...,f_p)^T: F\in \mathcal{N}_{\Psi}^{\otimes p}(\Omega), {\rm dim(span}\{f_1,...,f_p\}) = r\},
\end{align*}
and define
\begin{align*}
    \mathcal{B}_r(R) = \{f\in \mathcal{B}_r: \|F\|\leq R\}
\end{align*}
for $r=1,...,p$. Thus, $\mathcal{A}_r(R) = \cup_{k=1}^r \mathcal{B}_k(R)$ and $\mathcal{A}_r = \cup_{k=1}^r \mathcal{B}_k$. We start with some properties of $\mathcal{B}_r(R)$. The first proposition states that the orthogonal transformation does not change the norm, where the proof is provided in Appendix \ref{pfpropoth}.
\begin{proposition}\label{propothchange}
For any orthogonal matrix $U$ and $F\in \mathcal{B}_r(R)$, $UF\in \mathcal{B}_r(R)$ and $\|UF\|=\|F\|$.
\end{proposition}
Let $(T,d)$ be a metric space with metric $d$, and $T$ is a space. The $\epsilon$-covering number of the metric space $(T,d)$, denoted as $N(\epsilon,T,d)$, is the minimum integer $N$ so that there exist $N$ distinct balls in $(T,d)$ with radius $\epsilon$, and the union of these balls covers $T$. Let $H(\epsilon,T,d) = \log N(\epsilon,T,d)$ be the entropy number. Let $\|G\|_n$ be the emprical norm for $G\in \mathcal{B}_r(R)$, defined by
\begin{align*}
    \|G\|_n^2 = \frac{1}{n}\sum_{k=1}^n\sum_{j=1}^pg_j(x_k)^2,
\end{align*}
where $G = (g_1,...,g_p)^T$. In the next proposition, we provide an upper bound of the entropy number of  $\mathcal{B}_r(1)$, where the proof is provided in Appendix \ref{pfpropentB}.
\begin{proposition}\label{propentropy1}
The entropy number of $\mathcal{B}_r(1)$ is bounded by
\begin{align}\label{entropyb1}
     H(\delta, \mathcal{B}_r(1),\|\cdot\|_n) \leq & pr\log \left(1 + \frac{8\sqrt{r}}{\delta}\right) + C_1r\bigg(\frac{\sqrt{r}}{\delta}\bigg)^{d/\nu},
\end{align}
where $C_1$ is a constant.
\end{proposition}
Note $\mathcal{A}_r(1) = \cup_{k=1}^r \mathcal{B}_k(1)$. The following proposition is a direct result of Proposition \ref{propentropy1}.
\begin{proposition}\label{propentropyB}
The entropy number of $\mathcal{A}_r(1)$ is bounded by
\begin{align}\label{entropyb1}
     H(\delta, \mathcal{A}_r(1),\|\cdot\|_n) \leq & pr\log \left(1 + \frac{8\sqrt{r}}{\delta}\right) + Cr\bigg(\frac{\sqrt{r}}{\delta}\bigg)^{d/\nu},
\end{align}
where $C$ is a constant.
\end{proposition}

\section{Proof of Theorem \ref{mainthm1}}
Before we prove Theorem \ref{mainthm1} in the main text, we first present two lemmas. The first lemma is related to the empirical inner product between $f$ and $\epsilon$, and the second lemma states that the ratio of the empirical norm divided by the $L_2$ norm can be bounded.

We define the $L_\infty$ norm and $L_2$ norm for $G\in \mathcal{N}_{\Psi}^{\otimes p}(\Omega)$ as follows. Let $G=(g_1,...,g_p)^T \in \mathcal{N}_{\Psi}^{\otimes p}(\Omega)$. We define
\begin{align*}
    \|G\|_{L_\infty(\Omega)} = \sup_{1\leq j\leq p} \|g_j\|_{L_\infty(\Omega)},
\end{align*}
and
\begin{align*}
    \|G\|_{L_2(\Omega)}^2 = \sum_{j=1}^p \|g_j\|_{L_2(\Omega)}^2.
\end{align*}
\begin{lemma}\label{leminnersmall}
Suppose that $\epsilon_{kj}$'s are sub-Gaussian and $\sup_{G\in \mathcal{A}_r(1)}\|G\|_n\leq R$. Then we have
\begin{align*}
    P\left(\sup_{G\in \mathcal{A}_r(1)}\frac{\left|\frac{1}{pn}\sum_{k=1}^p\sum_{j=1}^n \epsilon_{kj}g_k(x_j)\right|}{\left(\sqrt{pr}\|G\|_n \log^{1/2}\left(1+\frac{4\sqrt{r}}{\|G\|_n}\right)+ r^{\frac{2\nu+d}{4\nu}} \|G\|_n^{1-\frac{d}{2\nu}}\right)} \geq T_1 \right) \leq c\exp(-npr^{\frac{2\nu+d}{2\nu}}T_1^2/c^2).
\end{align*}
for any $T_1 > 0$,  where $c>0$ is a constant.
\end{lemma}

\begin{lemma}\label{lemmaratio1}
Assume for class $\mathcal{G}\subset \mathcal{A}_r(1)$, $\sup_{G\in \mathcal{G}}\|G\|_{L_\infty(\Omega)}\leq K < 1$, and $n\delta_n\rightarrow \infty$. Then we have with probability at least $1-r\exp(-n\delta_n)$,
\begin{align*}
\sup_{G\in \mathcal{G}}\left|\|G\|^2_n-\|G\|^2_{L_2(\Omega)}\right| \leq Cr\delta_n.
\end{align*}
for some positive constant $C$ only depending on $\Omega$.
\end{lemma}

Now we are ready to prove Theorem \ref{mainthm1}. Let $r_2= r_1 + r$, where $r_1$ is as in \eqref{firstmodelXX} and $r = {\rm dim}(\mathcal{F})$. Let $\mathcal{F}_1 = {\rm span}\{ \hat f_1,...,\hat f_p\}$. Therefore, dim$(\mathcal{F}_1 + \mathcal{F})\leq r_2$, where
\begin{align*}
    \mathcal{F}_1 + \mathcal{F} = \{g|g = h_1 + h_2,h_1\in \mathcal{F}_1, h_2\in \mathcal{F}\}.
\end{align*}
Because $\hat F$ is the solution to \eqref{firstmodelXX}, we have
\begin{align*}
    \frac{1}{pn}\sum_{j=1}^p \|Y_j - \hat f_j(X)\|_2^2 + \lambda_1 \sum_{j=1}^p\|\hat f_j\|_{\mathcal{N}_{\Psi}(\Omega)}^2 \leq \frac{1}{pn}\sum_{j=1}^p \|Y_j - f_j(X)\|_2^2 + \lambda_1 \sum_{j=1}^p\|f_j\|_{\mathcal{N}_{\Psi}(\Omega)}^2,
\end{align*}
which is the same as
\begin{align}\label{consisineq}
     \frac{1}{pn}\sum_{j=1}^p \|f_j(X) - \hat f_j(X)\|_2^2 + \lambda_1 \sum_{j=1}^p\|\hat f_j\|_{\mathcal{N}_{\Psi}(\Omega)}^2 \leq  \frac{2}{pn}\sum_{j=1}^p \langle f_j(X) - \hat f_j(X),\epsilon_i\rangle_2 + \lambda_1 \sum_{j=1}^p\|f_j\|_{\mathcal{N}_{\Psi}(\Omega)}^2,
\end{align}
where $\langle a,b\rangle_2$ is the inner product of two vectors $a$ and $b$.

Applying Lemma \ref{leminnersmall} to the function class
\begin{align*}
    \left\{g|g = \frac{F - \hat F}{\|F\| + \|\hat F\|}, \hat F \mbox{ is the solution to \eqref{firstmodelXX}} \right\} \subset\mathcal{A}_{r_2}(1)
\end{align*}
yields
\begin{align}\label{consisXineq1}
    & \frac{1}{p}\|F - \hat F\|_n^2 + \lambda_1 \|\hat F\|^2\nonumber\\
    \leq & O_{P}((np)^{-1/2}r_2^{-\frac{2\nu+d}{4\nu}})\left(\sqrt{pr_2}\|F - \hat F\|_{n} \log^{1/2}\left(1+\frac{4\sqrt{r_2}}{\|F - \hat F\|_{n}}\right)+ r_2^{\frac{2\nu+d}{4\nu}} \|F - \hat F\|_{n}^{1-\frac{d}{2\nu}}(\|F\| + \|\hat F\|)^{\frac{d}{2\nu}}\right) \nonumber\\
    & + \lambda_1\|F\|^2.
\end{align}
We consider two cases.

\textbf{Case 1:} $\|F\| \leq \|\hat F\|$. By \eqref{consisXineq1}, we have
\begin{align}\label{consiscase11eq1}
    & \frac{1}{p}\|F - \hat F\|_n^2 + \lambda_1 \|\hat F\|^2\nonumber\\
    \leq & O_{P}((np)^{-1/2}r_2^{-\frac{2\nu+d}{4\nu}})\left(\sqrt{pr_2}\|F - \hat F\|_{n} \log^{1/2}\left(1+\frac{4\sqrt{r_2}}{\|F - \hat F\|_{n}}\right)+ r_2^{\frac{2\nu+d}{4\nu}} \|F - \hat F\|_{n}^{1-\frac{d}{2\nu}}\|\hat F\|^{\frac{d}{2\nu}}\right) + \lambda_1\|F\|^2.
\end{align}
Then we have either
\begin{align}\label{case11c1eq1}
    \frac{1}{p}\|F - \hat F\|_n^2 + \lambda_1 \|\hat F\|^2 \leq 4\lambda_1\|F\|^2
\end{align}
or
\begin{align}\label{case11eq2}
    & \frac{1}{p}\|F - \hat F\|_n^2 + \lambda_1 \|\hat F\|^2\nonumber\\
    \leq & O_{P}((np)^{-1/2}r_2^{-\frac{2\nu+d}{4\nu}})\left(\sqrt{pr_2}\|F - \hat F\|_{n} \log^{1/2}\left(1+\frac{4\sqrt{r_2}}{\|F - \hat F\|_{n}}\right)+ r_2^{\frac{2\nu+d}{4\nu}} \|F - \hat F\|_{n}^{1-\frac{d}{2\nu}}\|\hat F\|^{\frac{d}{2\nu}}\right).
\end{align}
It can be seen that \eqref{case11c1eq1} implies
\begin{align*}
    \|F - \hat F\|_n^2 \leq 4p\lambda_1\|F\|^2.
\end{align*}
Now consider \eqref{case11eq2}. We have two subcases.

\textbf{Case 1.1:}
\begin{align*}
    \sqrt{pr_2}\|F - \hat F\|_{n} \log^{1/2}\left(1+\frac{4\sqrt{r_2}}{\|F - \hat F\|_{n}}\right) \leq r_2^{\frac{2\nu+d}{4\nu}} \|F - \hat F\|_{n}^{1-\frac{d}{2\nu}}\|\hat F\|^{\frac{d}{2\nu}}.
\end{align*}
Then by \eqref{consiscase11eq1}, we have
\begin{align}\label{case111eq1}
    & \frac{1}{p}\|F - \hat F\|_n^2 + \lambda_1 \|\hat F\|^2\leq O_{P}((np)^{-1/2}) \|F - \hat F\|_{n}^{1-\frac{d}{2\nu}}\|\hat F\|^{\frac{d}{2\nu}}.
\end{align}
Solving \eqref{case111eq1}, we have
% \begin{align*}\label{case111eqz1}
%     \lambda_1 \|\hat F\|^2  \leq &  O_{P}((np)^{-1/2})r_2^{\frac{\nu+d}{2\nu}} \|F - \hat F\|_{n}^{1-\frac{d}{2\nu}}\|\hat F\|^{\frac{d}{2\nu}}\\
%     \Rightarrow & \|\hat F\|^{\frac{4\nu - d}{2\nu}}
% \end{align*}
\begin{align}\label{case111eqre}
    \|F - \hat F\|_n^2 = O_{P}(n^{-1}p^{\frac{2\nu-d}{2\nu}}\lambda_1^{-\frac{d}{2\nu}}).
\end{align}

\textbf{Case 1.2:}
\begin{align*}
    \sqrt{pr_2}\|F - \hat F\|_{n} \log^{1/2}\left(1+\frac{4\sqrt{r_2}}{\|F - \hat F\|_{n}}\right) > r_2^{\frac{2\nu+d}{4\nu}} \|F - \hat F\|_{n}^{1-\frac{d}{2\nu}}\|\hat F\|^{\frac{d}{2\nu}}.
\end{align*}
This implies
\begin{align}\label{case12eq1}
    \sqrt{p}\|F - \hat F\|_{n}^{\frac{d}{2\nu}} \log^{1/2}\left(1+\frac{4\sqrt{r_2}}{\|F - \hat F\|_{n}}\right) > r_2^{\frac{d}{4\nu}} \|F\|^{\frac{d}{2\nu}}.
\end{align}
Under Case 1.2, by \eqref{case11eq2}, we have
\begin{align*}
    & \frac{1}{p}\|F - \hat F\|_n^2 + \lambda_1 \|\hat F\|^2\leq  O_{P}(n^{-1/2}r_2^{-\frac{2\nu+d}{4\nu}})\sqrt{r_2}\|F - \hat F\|_{n} \log^{1/2}\left(1+\frac{4\sqrt{r_2}}{\|F - \hat F\|_{n}}\right),
\end{align*}
which leads to
\begin{align}\label{case12eq2}
    & \|F - \hat F\|_n\leq  O_{P}(pn^{-1/2}r_2^{-\frac{d}{4\nu}}) \log^{1/2}\left(1+\frac{4\sqrt{r_2}}{\|F - \hat F\|_{n}}\right).
\end{align}
Combining \eqref{case12eq1} and \eqref{case12eq2}, we have
\begin{align}\label{case12eq3}
    r_2^{\frac{d}{4\nu}} \|F\|^{\frac{d}{2\nu}} < O_{P}(n^{-\frac{d}{4\nu}}p^{\frac{d}{2\nu} +\frac{1}{2}}r_2^{-\frac{d^2}{8\nu^2}}) \log^{\frac{d}{4\nu}+\frac{1}{2}}\left(1+\frac{4\sqrt{r_2}}{\|F - \hat F\|_{n}}\right) .
\end{align}
Solving \eqref{case12eq3}, we have for some constant $C$
\begin{align}\label{case12re}
    \|F - \hat F\|_{n} < 4\sqrt{r_2}\exp\left[-Cr_2^{\frac{d}{2\nu}} \|F\|^{\frac{2d}{2\nu+d}}n^{\frac{d}{2\nu+d}}p^{-\frac{2d+2\nu}{2\nu+d}} \right],
\end{align}
which contradicts \eqref{case12eq1} because $r_1^{\frac{d}{2\nu}} \|F\|^{\frac{2d}{2\nu+d}}n^{\frac{d}{2\nu+d}}p^{-\frac{2d+2\nu}{2\nu+d}}$ goes to infinity and $r\leq r_1$.

\textbf{Case 2:} $\|F\| > \|\hat F\|$. By \eqref{consisXineq1}, we have
\begin{align}\label{consiscase2eq1}
    & \frac{1}{p}\|F - \hat F\|_n^2 + \lambda_1 \|\hat F\|^2\nonumber\\
    \leq & O_{P}((np)^{-1/2}r_2^{-\frac{2\nu+d}{4\nu}})\left(\sqrt{pr_2}\|F - \hat F\|_{n} \log^{1/2}\left(1+\frac{4\sqrt{r_2}}{\|F - \hat F\|_{n}}\right)+ r_2^{\frac{2\nu+d}{4\nu}} \|F - \hat F\|_{n}^{1-\frac{d}{2\nu}}\|F\|^{\frac{d}{2\nu}}\right)\nonumber\\
    & + \lambda_1\|F\|^2.
\end{align}
Then we have either
\begin{align}\label{case2c1eq1}
    \frac{1}{p}\|F - \hat F\|_n^2 + \lambda_1 \|\hat F\|^2 \leq 4\lambda_1\|F\|^2
\end{align}
or
\begin{align}\label{case2c2eq2}
    & \frac{1}{p}\|F - \hat F\|_n^2 + \lambda_1 \|\hat F\|^2\nonumber\\
    \leq & O_{P}((np)^{-1/2}r_2^{-\frac{2\nu+d}{4\nu}})\left(\sqrt{pr_2}\|F - \hat F\|_{n} \log^{1/2}\left(1+\frac{4\sqrt{r_2}}{\|F - \hat F\|_{n}}\right)+ r_2^{\frac{2\nu+d}{4\nu}} \|F - \hat F\|_{n}^{1-\frac{d}{2\nu}}\|\hat F\|^{\frac{d}{2\nu}}\right).
\end{align}
It can be seen that \eqref{case2c1eq1} implies
\begin{align*}
    \|F - \hat F\|_n^2 \leq 4p\lambda_1\|F\|^2.
\end{align*}
Now consider \eqref{case2c2eq2}. We have two subcases.

\textbf{Case 2.1:}
\begin{align*}
    \sqrt{pr_2}\|F - \hat F\|_{n} \log^{1/2}\left(1+\frac{4\sqrt{r_2}}{\|F - \hat F\|_{n}}\right) \leq r_2^{\frac{2\nu+d}{4\nu}} \|F - \hat F\|_{n}^{1-\frac{d}{2\nu}}\|F\|^{\frac{d}{2\nu}}.
\end{align*}
Then by \eqref{consiscase2eq1}, we have
\begin{align}\label{case21eq1}
    & \frac{1}{p}\|F - \hat F\|_n^2 + \lambda_1 \|\hat F\|^2\leq O_{P}((np)^{-1/2}) \|F - \hat F\|_{n}^{1-\frac{d}{2\nu}}\|F\|^{\frac{d}{2\nu}}.
\end{align}
Solving \eqref{case21eq1}, we have
\begin{align}\label{case2eqre}
    \|F - \hat F\|_n^2 = O_{P}\left(\left(\frac{p}{n}\right)^{\frac{2\nu}{2\nu+d}}\right)\|F\|^{\frac{2d}{2\nu+d}}.
\end{align}

\textbf{Case 2.2:}
The case
\begin{align*}
    \sqrt{pr_2}\|F - \hat F\|_{n} \log^{1/2}\left(1+\frac{4\sqrt{r_2}}{\|F - \hat F\|_{n}}\right) > r_2^{\frac{2\nu+d}{4\nu}} \|F - \hat F\|_{n}^{1-\frac{d}{2\nu}}\|F\|^{\frac{d}{2\nu}}
\end{align*}
is similar to Case 1.2.

By choosing  $$\lambda_1 \asymp n^{-\frac{2\nu}{2\nu + d}}p^{-\frac{d}{2\nu + d}} \|F\|^{-\frac{4\nu}{2\nu+d}},$$ we obtain
\begin{align*}
    \|F - \hat F\|_n^2 =O_{P}\left(\left(\frac{p}{n}\right)^{\frac{2\nu}{2\nu+d}}\right)\|F\|^{\frac{2d}{2\nu+d}}.
\end{align*}
For some constant $C_1$ depending on $\Omega$ and $\Psi$, let $\delta_n \geq C_1n^{-\nu/(2\nu + d)}$. Noting that $r_1\geq r$, and choosing $\delta_n = n^{-\frac{2\nu}{2\nu+d}}$ in Lemma \ref{lemmaratio1}, we have
\begin{align*}
    \|F - \hat F\|_{L_2(\Omega)}^2 \leq \|F - \hat F\|_{n}^2 + C_2r_1n^{-\frac{2\nu}{2\nu+d}} = O_P(Sn^{-\frac{2\nu}{2\nu+d}}),
\end{align*}
where $S$ is as in Theorem \ref{mainthm1}. Thus, we obtain the $L_2$ convergence rate as desired.

\section{Proofs of Propositions and Lemmas}

\subsection{Proof of Proposition \ref{propothchange}}\label{pfpropoth}
Let $F=(f_1,...,f_p)^T\in \mathcal{B}_r(R)$. Because any orthogonal transformation does not change the rank, it can be seen that dim$($span$(UF))= r$. Thus, it is enough to show $\|UF\|=\|F\|$, because $\|F\|\leq R$. The norm of $\|UF\|$ can be calculated by
\begin{align*}
    \|UF\|^2 = & \sum_{j=1}^p \|\sum_{k=1}^p u_{jk}f_k \|_{\mathcal{N}_{\Psi}(\Omega)}^2
    =  \sum_{j=1}^p \langle\sum_{k=1}^p u_{jk}f_k,\sum_{l=1}^p u_{jl}f_l \rangle_{\mathcal{N}_{\Psi}(\Omega)}\\
    = & \sum_{j=1}^p \sum_{k,l=1}^p u_{jk}u_{jl}\langle f_k,f_l\rangle_{\mathcal{N}_{\Psi}(\Omega)}
    =  \sum_{k,l=1}^p \sum_{j=1}^p  u_{jk}u_{jl}\langle f_k,f_l\rangle_{\mathcal{N}_{\Psi}(\Omega)}\\
    = & \sum_{k=1}^p \|f_k\|_{\mathcal{N}_{\Psi}(\Omega)}^2=\|F\|,
\end{align*}
where the fifth equality is because $U$ is orthogonal.

\subsection{Proof of Proposition \ref{propentropy1}}\label{pfpropentB}
For $F\in \mathcal{B}_r(R)$, there exists $r$ linear independent basis. We can use a permutation matrix $P\in \RR^{r\times p}$ such that elements in $PF$ are linear independent. Then we can write $F$ as
\begin{align*}
    F = APF,
\end{align*}
where $A\in \RR^{p\times r}$, because any $f_j\in F$ can be written as a linear combination of $PF$. Now we apply singular value decomposition to $A$ such that $A = U\Lambda V^T$, where
\begin{align*}
    \Lambda = \left(\begin{array}{cc}
         \Lambda_1  \\
         0
    \end{array}\right)\in \RR^{p\times r},
\end{align*}
with $\Lambda_1 = {\rm diag}(\lambda_1,...,\lambda_r)$ are singular values, and $U\in \RR^{p\times p}$ and $V \in \RR^{r\times r}$ are two orthogonal matrices. Putting all things together, we have
\begin{align*}
    F = U\Lambda V^TPF = U V_r,
\end{align*}
where $V_r = \Lambda V^TPF = (\lambda_1v_1,...,\lambda_rv_r,0,...,0)^T$. Let $w_1 = (\lambda_1v_1,...,\lambda_rv_r)^T$. By Proposition \ref{propothchange}, we have $\|F\| = \|V_r\|$. Thus, $\|V_r\|\leq R$. Write
\begin{align*}
    U = \left(\begin{array}{cc}
        u_{11}^T & u_{12}^T \\
        u_{21}^T & u_{22}^T\\
        \vdots & \vdots\\
        u_{p1}^T & u_{p2}^T
    \end{array}\right),
\end{align*}
where $u_{j1} \in \RR^r$ and $u_{j2} \in \RR^{p-r}$ for $j=1,...,p$. Thus,
\begin{align*}
    F = (u_{11}^T w_1,...,u_{p1}^T w_1)^T.
\end{align*}
Let $\mathcal{N}_{\Psi}(R) = \{f\in \mathcal{N}_{\Psi}(\Omega);\|f\|_{\mathcal{N}_{\Psi}(\Omega)}\leq 1\}$. Since $\mathcal{N}_{\Psi}(\Omega)$ is equivalent to the Sobolev space $H^{\nu}(\Omega)$ (\cite{wendland2004scattered}, Corollary 10.13), the entropy number of $\mathcal{N}_{\Psi}(R)$ can be bounded by \citep{adams2003sobolev}
\begin{align*}
H(\delta,\mathcal{N}_{\Psi}(1),\|\cdot\|_{L_\infty(\Omega)})\leq C\bigg(\frac{R_1}{\delta}\bigg)^{d/\nu},
\end{align*}
where $C$ and $R_1$ are two positive constants. Thus, the entropy number of set
\begin{align*}
    \mathcal{W} = \{w = (w_1,...,w_r)^T: \sum_{j=1}^r \|w_j\|_{\mathcal{N}_{\Psi}(\Omega)}^2 \leq 1\}
\end{align*}
can be bounded by
\begin{align*}
H(\delta,\mathcal{W},\|\cdot\|_\infty)\leq Cr\bigg(\frac{R_1}{\delta}\bigg)^{d/\nu},
\end{align*}
where $\|\cdot\|_\infty$ for function class $\mathcal{W}$ is defined by
\begin{align*}
    \|w\|_\infty = \max_{1\leq j\leq r} \|w_j\|_{L_\infty(\Omega)},
\end{align*}
for $w\in \mathcal{W}$.

Next, consider matrix $U_r = (u_{11},...,u_{p1})^T \in \RR_2^{p \times r}$. Because $U$ is orthogonal, $U_r^TU_r = I_r$, where $I_r$ is the identity matrix. Therefore, $\|U_r\|_\infty\leq \|U_r\|_F \leq \sqrt{r}\|U_r\|_2 = \sqrt{r}$, where $\|U_r\|_F$ is the Frobenius norm of $U_r$.

Now consider the covering number of $\mathcal{B}_r(1)$. For $F,F'\in \mathcal{B}_r(1)$, we can write $F = U_rw_1$ and $F' = U_r'w_2'$, where $U_r'\in \RR_2^{p\times r}$ and $w_2'\in \mathcal{W}$. Thus, by the triangel inequality, we have
\begin{align}\label{ineqforentropy1}
    \|F - F'\|_n & = \|U_rw_1 - U_r'w_2'\|_n = \|U_rw_1 -U_r'w_1+U_r'w_1- U_r'w_2'\|_n\nonumber\\
    & \leq \|U_rw_1 -U_r'w_1\|_n+\|U_r'w_1 - U_r'w_2'\|_n \leq \|U_r -U_r'\|_2\|w_1\|_n + \|w_1 -w_2\|_n\nonumber \\
    & \leq \|U_r -U_r'\|_F\|w_1\|_n + \sqrt{r}\|w_1 -w_2\|_\infty  \leq\|U_r -U_r'\|_F +  \sqrt{r}\|w_1 - w_2\|_\infty.
\end{align}
Note that the entropy number of set
\begin{align*}
    \mathcal{U} = \{U_r\in \RR^{p\times r}: \|U_r\|_F\leq \sqrt{r}\}
\end{align*}
can be bounded by
\begin{align*}
    H(\delta,\mathcal{U},\|\cdot\|_F)\leq pr\log \left(1 + \frac{4\sqrt{r_2}}{\delta}\right).
\end{align*}
Therefore, by \eqref{ineqforentropy1}, the entropy number of $\mathcal{B}_r(1)$ can be bounded by
\begin{align}\label{entropy1}
    H(\delta, \mathcal{B}_r(1),\|\cdot\|_n)\leq &  H(\delta/2,\mathcal{U},\|\cdot\|_F) + H(\delta/(2\sqrt{r}),\mathcal{W},\|\cdot\|)\nonumber\\
    \leq & pr\log \left(1 + \frac{8\sqrt{r}}{\delta}\right) + Cr\bigg(\frac{\sqrt{r}R_1}{\delta}\bigg)^{d/\nu},
\end{align}
which finishes the proof.

\subsection{Proof of Lemma \ref{leminnersmall}}
Before the proof, we first present a lemma used in this proof.
\begin{lemma}[Corollary 8.3 of \cite{geer2000empirical}]\label{coro83invan}
Suppose that $\epsilon_j$'s are sub-Gaussian and $\sup_{g\in \mathcal{G}}\|g\|_n\leq R$, where
\begin{align*}
    \|g\|_n^2 = \frac{1}{n}\sum_{j=1}^ng(x_j)^2.
\end{align*}
Suppose
\begin{align*}
    \int_0^R H^{1/2}(u,\mathcal{G},\|\cdot\|_n)du < \infty.
\end{align*}
Then for some constant $C$ depending only on the parameters of sub-Gaussian random variables, and for $\delta>0$ and
\begin{align}
    \sqrt{n}\delta \geq 2C\max\left\{\int_{0}^R H^{1/2}(u,\mathcal{G},\|\cdot\|_n)du, R\right\},
\end{align}
we have
\begin{align}
    P\left(\sup_{g\in \mathcal{G}}\left|\frac{1}{n}\sum_{j=1}^n \epsilon_jg(x_j)\right|\geq \delta \right)\leq C\exp\left[-\frac{n\delta^2}{4C^2R^2}\right].
\end{align}
\end{lemma}

Now we are ready to prove Lemma \ref{leminnersmall}. By Proposition \ref{propentropyB}, we have
\begin{align*}
    H(\delta,\mathcal{A}_r(1),\|\cdot\|_n)\leq  pr\log \left(1 + \frac{8r}{\delta}\right) + Cr\bigg(\frac{\sqrt{r}}{\delta}\bigg)^{d/\nu}.
\end{align*}

Therefore, we have
\begin{align}\label{pfpeel1ent}
    \int_{0}^\delta H^{1/2}(u, \mathcal{A}_r(1),\|\cdot\|_n)du \leq & \int_{0}^\delta \sqrt{pr\log \left(1 + \frac{4\sqrt{r}}{u}\right) + Cr\bigg(\frac{\sqrt{r}}{u}\bigg)^{d/\nu}}du\nonumber\\
    \leq & \int_{0}^\delta \sqrt{pr\log \left(1 + \frac{4\sqrt{r}}{u}\right)}du + \sqrt{Cr}\bigg(\frac{\sqrt{r}}{u}\bigg)^{d/(2\nu)}du\nonumber\\
    \leq & \delta^{1/2}\left(\int_{0}^\delta pr\log \left(1 + \frac{4\sqrt{r}}{u}\right)du\right)^{1/2} + C_1r^{\frac{2\nu+d}{4\nu}} \delta^{1-\frac{d}{2\nu}}\nonumber\\
    = & \sqrt{pr}\delta^{1/2}\left(\delta \log\left(1+\frac{4\sqrt{r}}{\delta}\right) +  4\sqrt{r}\log\left(1 + \frac{\delta}{4\sqrt{r}}\right)  \right)^{1/2} + C_1r^{\frac{2\nu+d}{4\nu}} \delta^{1-\frac{d}{2\nu}}\nonumber\\
    \leq & \sqrt{pr}\delta \left( \log\left(1+\frac{4\sqrt{r}}{\delta}\right) +  1\right)^{1/2} + C_1r^{\frac{2\nu+d}{4\nu}} \delta^{1-\frac{d}{2\nu}}\nonumber\\
    \leq & C_2\sqrt{pr}\delta \log^{1/2}\left(1+\frac{4\sqrt{r}}{\delta}\right) + C_1r^{\frac{2\nu+d}{4\nu}} \delta^{1-\frac{d}{2\nu}}.
\end{align}
In \eqref{pfpeel1ent}, the first inequality is by Proposition \ref{propentropyB}; the second inequality is by the basic inequality $\sqrt{a+b} \leq \sqrt{a} + \sqrt{b}$ for $a,b>0$; the third inequality is by the Cauchy-Schwarz inequality; the fourth inequality is by $\log(1+x)\leq x$ for $x>0$.

Thus, for $T\geq 1$, by Lemma \ref{coro83invan}, let $C_3 = C_1\vee C_2$, we have
\begin{align*}
     & P\left(\sup_{g\in \mathcal{A}_r(1), \|g\|_n\leq \delta} \left|\frac{1}{pn}\sum_{k=1}^p\sum_{j=1}^n \epsilon_{kj}g_k(x_j)\right|\geq 2CTC_3\left(\sqrt{pr}\delta \log^{1/2}\left(1+\frac{4\sqrt{r}}{\delta}\right) + r^{\frac{2\nu+d}{4\nu}} \delta^{1-\frac{d}{2\nu}}\right) \right) \\
    \leq & C\exp\left[-npC_3^2T^2 \left(pr \log\left(1+\frac{4\sqrt{r}}{\delta}\right) + r^{\frac{2\nu+d}{2\nu}} \delta^{-\frac{d}{\nu}}\right)\right].
\end{align*}
Let $\mathcal{Q}_s = \{g\in \mathcal{A}_r(1), 2^{-s}R \leq \|g\|_n \leq 2^{-s+1}R\}$, for $s=1,2,...$. Therefore, we have
\begin{align*}
    & P\left(\sup_{g\in \mathcal{A}_r(1)}\frac{\left|\frac{1}{pn}\sum_{k=1}^p\sum_{j=1}^n \epsilon_{kj}g_k(x_j)\right|}{\left(\sqrt{pr}\|g\|_n \log^{1/2}\left(1+\frac{4\sqrt{r}}{\|g\|_n}\right)+ r^{\frac{2\nu+d}{4\nu}} \|g\|_n^{1-\frac{d}{2\nu}}\right)} \geq T_1 \right)\\
    \leq & \sum_{s=1}^\infty P\left(\sup_{g\in \mathcal{A}_r(1),g\in \mathcal{Q}_s}\frac{\left|\frac{1}{pn}\sum_{k=1}^p\sum_{j=1}^n \epsilon_{kj}g_k(x_j)\right|}{\left(\sqrt{pr}\|g\|_n \log^{1/2}\left(1+\frac{4\sqrt{r}}{\|g\|_n}\right)+ r^{\frac{2\nu+d}{4\nu}} \|g\|_n^{1-\frac{d}{2\nu}}\right)} \geq T_1 \right)\\
    \leq & \sum_{s=1}^\infty P\left(\sup_{g\in \mathcal{A}_r(1),g\in \mathcal{Q}_s}\left|\frac{1}{pn}\sum_{k=1}^p\sum_{j=1}^n \epsilon_{kj}g_k(x_j)\right|\geq T_1 \left(\sqrt{pr}2^{-s}R \log^{1/2}\left(1+\frac{4\sqrt{r}}{2^{-s+1}R}\right)+ r^{\frac{2\nu+d}{4\nu}} (2^{-s}R)^{1-\frac{d}{2\nu}}\right)  \right)\\
    \leq & \sum_{s=1}^\infty P\bigg(\sup_{g\in \mathcal{A}_r(1),g\in \mathcal{Q}_s}\left|\frac{1}{pn}\sum_{k=1}^p\sum_{j=1}^n \epsilon_{kj}g_k(x_j)\right|\geq 2^{1-\frac{d}{2\nu}}T_1 \bigg(\sqrt{pr}2^{-s+1}R \log^{1/2}\left(1+\frac{4\sqrt{r}}{2^{-s+1}R}\right)\\
    & + r^{\frac{2\nu+d}{4\nu}} (2^{-s+1}R)^{1-\frac{d}{2\nu}}\bigg)  \bigg)\\
    \leq & \sum_{s=1}^\infty C\exp\left[-np 2^{-d/\nu} T_1^2 C^{-2} \left(pr \log\left(1+\frac{4\sqrt{r}}{2^{-s+1}R}\right) + r^{\frac{2\nu+d}{2\nu}} (2^{-s+1}R)^{-\frac{d}{\nu}}\right)\right]\\
    \leq & C_4 \exp(-npr^{\frac{2\nu+d}{2\nu}}T_1^2/C_4^2),
\end{align*}
where the last inequality is because $p\geq r$. Therefore, we finish the proof.

\subsection{Proof of Lemma \ref{lemmaratio1}}
We need the following lemma, which is a direct result of Theorem 2.1 of \cite{van2014uniform}. Lemma \ref{thm21inGeer2014} provides an upper bound on the difference between the empirical norm and $L_2$ norm. In Lemma \ref{thm21inGeer2014}, we use the following definition. For $z>0$, we define
\begin{align*}
    J_\infty^2(z,\mathcal{G}_0) = C_0^2\inf_{\delta>0} \mathbb{E}\left[z\int_\delta^1 \sqrt{H(uz/2,\mathcal{G}_0,\|\cdot\|_\infty)}du + \sqrt{n}\delta z \right]^2,
\end{align*}
where $C_0$ is a constant, and $H(u,\mathcal{G}_0,\|\cdot\|_\infty)$ is the entropy of $(\mathcal{G}_0,\|\cdot\|_\infty)$ for a function class $\mathcal{G}_0$.
\begin{lemma}\label{thm21inGeer2014}
Let $\mathcal{G}_0 \subset \{g| \|g\|_{\mathcal{N}_\Psi(\Omega)} \leq 1\}$,  $R=\sup_{f\in\mathcal{G}_0}\|f\|_{L_2(\Omega)}$, and $K=\sup_{f\in\mathcal{G}_0}\|f\|_\infty$. Then for all $t>0$, with probability at least $1-\exp(-t)$,
\begin{align*}
\sup_{f\in\mathcal{G}_0}\bigg|\|f\|^2_n-\|f\|^2_{L_2(\Omega)}\bigg|\leq C_1\bigg(\frac{2R J_\infty(K,\mathcal{G}_0)+RK\sqrt{t}}{\sqrt{n}}+\frac{4 J_\infty^2(K,\mathcal{G}_0)+K^2t}{n}\bigg),
\end{align*}
where $C_1$ is a constant.
\end{lemma}
Now we are ready to prove Lemma \ref{lemmaratio1}. Similar to the proof of Proposition \ref{propothchange}, we can show that for any orthogonal transformation $U$, $\|UG\|_{L_2(\Omega)}^2 = \|G\|_{L_2(\Omega)}^2$. Then by the proof of Proposition \ref{propentropyB}, we can write $G$ as
\begin{align*}
    G = UV_r,
\end{align*}
where $V_r = (\lambda_1v_1,...,\lambda_rv_r,0,...,0)^T$ and $U$ is an orthogonal transformation. Thus, we have $\|G\|_{L_2(\Omega)}^2 = \|V_r\|_{L_2(\Omega)}^2$. Similarly, we have $\|G\|_n^2= \|V_r\|_n^2$. For each $\lambda_jv_j\in V_r$, by the interpolation inequality, we have
\begin{align*}
    \|\lambda_jv_j\|_{L_2(\Omega)}^2 \leq \|V_r\|_{L_2(\Omega)}^2 = \|G\|_{L_2(\Omega)}^2\leq \|G\|^2 \leq 1,
\end{align*}
and
\begin{align*}
    \|\lambda_jv_j\|_{\mathcal{N}_{\Psi}(\Omega)}^2 \leq \|V_r\|^2 \leq \|G\|^2 \leq 1.
\end{align*}
Since $\mathcal{N}_{\Psi}(\Omega)$ is equivalent to the Sobolev space $H^{\nu}(\Omega)$ (\cite{wendland2004scattered}, Corollary 10.13), the entropy number of $\mathcal{G}_j$ can be bounded by \citep{adams2003sobolev}
\begin{align*}
H(\delta,\mathcal{G}_j,\|\cdot\|_{L_\infty(\Omega)})\leq C\bigg(\frac{R_1}{\delta}\bigg)^{d/\nu}.
\end{align*}
The quantity $J_\infty(K,\mathcal{G}_j)$ then can be bounded by
\begin{align*}
    J_\infty^2(K,\mathcal{G}_0) \leq & C_0^2\left[K\int_0^1 \sqrt{H(uK/2,\mathcal{G}_0,\|\cdot\|_\infty)}du \right]^2\\
    \leq &  C_0^2\left[K\int_0^1 C\bigg(\frac{2R_1}{uK}\bigg)^{\frac{d}{2\nu}}du \right]^2\\
    \leq & C_1 K^{2-\frac{d}{\nu}}.
\end{align*}
Take $K=1$ and $R=1$. By Lemma \ref{thm21inGeer2014}, for a single class $\mathcal{G}_j$ that $\lambda_jv_j$ lies in, we have with probability at least $1-\exp(-t)$,
\begin{align*}
\sup_{f\in\mathcal{G}_j}\bigg|\|f\|^2_n-\|f\|^2_{L_2(\Omega)}\bigg|\leq C_2\bigg(\frac{2 +\sqrt{t}}{\sqrt{n}}+\frac{4+t}{n}\bigg).
\end{align*}
Therefore, by the union bound, and taking $t=n\delta_n$, we have with probability at least $1-r\exp(-n\delta_n)$,
\begin{align*}
    \sup_{G\in \mathcal{G}}\left|\|G\|^2_n-\|G\|^2_{L_2(\Omega)}\right| = & \sup_{G\in \mathcal{G}}\left|\|V_r\|^2_n-\|V_r\|^2_{L_2(\Omega)}\right| \leq \sum_{j=1}^r \sup_{f\in\mathcal{G}_j}\bigg|\|f\|^2_n-\|f\|^2_{L_2(\Omega)}\bigg|\nonumber\\
    \leq & C_2r\bigg(\frac{2 +\sqrt{t}}{\sqrt{n}}+\frac{4+t}{n}\bigg) = C_3r\delta_n.
\end{align*}
Thus, we finish the proof.

\bibliography{ref}
\end{document}